\ifx\shlhetal\undefinedcontrolsequence\let\shlhetal\relax\fi

\documentclass[11pt]{amsart}

\usepackage{amsmath}
\usepackage{amssymb}

\newtheorem{theorem}{Theorem}[section]
\newtheorem{claim}[theorem]{Claim}

\theoremstyle{definition}
\newtheorem{definition}[theorem]{Definition}

\theoremstyle{remark}
\newtheorem{remark}[theorem]{Remark}

\newcount\skewfactor
\def\mathunderaccent#1#2 {\let\theaccent#1\skewfactor#2
\mathpalette\putaccentunder}
\def\putaccentunder#1#2{\oalign{$#1#2$\crcr\hidewidth
\vbox to.2ex{\hbox{$#1\skew\skewfactor\theaccent{}$}\vss}\hidewidth}}
\def\name{\mathunderaccent\tilde-3 }

\def\smallbox#1{\leavevmode\thinspace\hbox{\vrule\vtop{\vbox
   {\hrule\kern1pt\hbox{\vphantom{\tt/}\thinspace{\tt#1}\thinspace}}
   \kern1pt\hrule}\vrule}\thinspace}

\newcommand{\cf}{{\rm cf}}

\newcommand{\Then}{{\underline{Then}}}

\def\qedref#1{$\qed_{\reforiginal{#1}}$}

\title{Strong polarized relations for the continuum}
\author{Shimon Garti}
\address{Institute of Mathematics
 The Hebrew University of Jerusalem
 Jerusalem 91904, Israel}
\email{shimon.garty@mail.huji.ac.il}

\author{Saharon Shelah}
\address{Institute of Mathematics
 The Hebrew University of Jerusalem
 Jerusalem 91904, Israel
 and  Department of Mathematics
 Rutgers University
 New Brunswick, NJ 08854, USA}
\email{shelah@math.huji.ac.il}
\urladdr{http://www.math.rutgers.edu/\char`\~shelah}

\thanks{First typed: February 2010 \newline Research supported by the United States-Israel Binational Science Foundation. Publication 964 of the second author}

\begin{document}
\let\labeloriginal\label
\let\reforiginal\ref

\begin{abstract}
We prove that the strong polarized relation $\binom{2^\mu}{\mu} \rightarrow \binom{2^\mu}{\mu}^{1,1}_2$ is consistent with ZFC. We show this for $\mu = \aleph_0$, and for every supercompact cardinal $\mu$. We also characterize the polarized relation below the splitting number.
\end{abstract}

\maketitle

\newpage

\section{introduction}

The balanced polarized relation $\binom{\alpha}{\beta} \rightarrow \binom{\gamma}{\delta}^{1,1}_2$ asserts that for every coloring $c : \alpha \times \beta \rightarrow 2$ there are $A \subseteq \alpha$ and $B \subseteq \beta$ such that ${\rm otp} (A) = \gamma, {\rm otp} (B) = \delta$ and $c \upharpoonright (A \times B)$ is constant. This relation was first introduced in \cite{MR0081864}, and investigated further in \cite{MR0202613}. A wonderful summary of the basic facts for this relation, appears in \cite{williams}. \newline 
Apparently, this relation can be true only when $\gamma \leq \alpha$ and $\delta \leq \beta$. It means that the strongest form of it is the case of $\binom{\alpha}{\beta} \rightarrow \binom{\alpha}{\beta}^{1,1}_2$. We can give a name to this situation:

\begin{definition}
\label{stpolarized}
The strong polarized relation. \newline 
If $\binom{\alpha}{\beta} \rightarrow \binom{\alpha}{\beta}^{1,1}_\theta$ then we say that the pair $(\alpha, \beta)$ satisfies the strong polarized relation with $\theta$ colors.
\end{definition}

If $2^\kappa = \kappa^+$ then $\binom{\kappa^+}{\kappa} \nrightarrow \binom{\kappa^+}{\kappa}^{1,1}_2$. This result, and similar negative results, go back to Sierpinsky. Despite the negative results under the (local) assumption of the GCH, we can show that the positive relation $\binom{\kappa^+}{\kappa} \rightarrow \binom{\kappa^+}{\kappa}^{1,1}_2$ is consistent with ZFC. Such a result appears in \cite{GaSh949}, and it was known for the specific case of $\kappa = \aleph_0$ (under the appropriate assumption, e.g., MA + $2^{\aleph_0} > \aleph_1$, see Laver in \cite{MR0371652} which proves that if Martin's Axiom holds for $\kappa$ then $\binom{\kappa}{\omega} \rightarrow \binom{\kappa}{\omega}^{1,1}_2$).

So the negative result under the GCH cannot become a theorem in ZFC. But we can restate this result in the form $\binom{2^\kappa}{\kappa} \nrightarrow \binom{2^\kappa}{\kappa}^{1,1}_2$. In this light, a natural question is whether the positive relation $\binom{2^\kappa}{\kappa} \rightarrow \binom{2^\kappa}{\kappa}^{1,1}_2$ is consistent with ZFC. One might suspect that the answer is negative, and this is the correct generalization of the negative result under the GCH. Notice that the result of Laver (from \cite{MR0371652}) does not help in this case, since ${\rm MA}_{2^{\aleph_0}}$ never holds.

Moreover, in the first section we prove that if $\aleph_0 < \cf(\mu) \leq \mu < \mathfrak{s}$, then $\binom{\mu}{\omega} \rightarrow \binom{\mu}{\omega}^{1,1}_2$. Hence, for every $\mu$ below the continuum (with uncountable cofinality), we can force a positive result by increasing $\mathfrak{s}$. But $\mathfrak{s} \leq 2^{\aleph_0}$, so this method does not help in our question. Again, under some assumptions one can prove negative results. For example, if $2^\omega$ is regular and MA(countable) holds, then $\binom{2^\omega}{\omega} \nrightarrow \binom{2^\omega}{\omega}^{1,1}_2$.

Nevertheless, we shall prove that a positive relation is consistent here. In the first section we deal with $\mu = \aleph_0$. Here we can use a finite support iteration of $ccc$ forcing notions, yielding a $ccc$ forcing notion in the limit of the sequence. We indicate that $\cf(2^\omega) \geq \aleph_2$ in our construction, and it might be that the relation $\binom{2^\omega}{\omega} \nrightarrow \binom{2^\omega}{\omega}^{1,1}_2$ is provable in ZFC whenever $\cf(2^\omega) = \aleph_1$. Notice that $\cf(2^\omega) = \aleph_1$ implies the existence of weak diamond on $\aleph_1$, so the negative relation becomes plausible.

In the second section we try to generalize it to higher cardinals. Here we encounter some difficulty, since the chain condition is not inherited from the members of the iteration. Starting with Laver-indestructible supercompact cardinal, we can overcome this problem (as well as other obstacles in the generalization of the countable case).

We try to use standard notation. We use the letters $\theta, \kappa, \lambda, \mu, \chi$ for infinite cardinals, and $\alpha, \beta, \gamma, \delta, \varepsilon, \zeta$ for ordinals. For a regular cardinal $\kappa$ we denote the ideal of bounded subsets of $\kappa$ by $J^{\rm bd}_\kappa$. For $A,B \subseteq \kappa$ we say that $A \subseteq^* B$ when $A \setminus B$ is bounded in $\kappa$. The symbol $[\lambda]^\kappa$ means the collection of all the subsets of $\lambda$ of cardinality $\kappa$. We denote the continuum by $\mathfrak{c}$.

Recall that by Laver (in \cite{MR0472529}) we can make a supercompact cardinal $\mu$ indestructible, upon forcing with $\mu$-closed forcing notions. We shall use this assumption in \S 2. We indicate that $p \leq q$ means (in this paper) that $q$ gives more information than $p$ in forcing notions.

An important specific forcing to be mentioned here is the Mathias forcing. Let $D$ be a nonprincipal ultrafilter on $\omega$. We define $\mathbb{M}_D$ as follows. The conditions in $\mathbb{M}_D$ are pairs of the form $(s, A)$ when $s \in [\omega]^{< \omega}$ and $A \in D$. For the order, $(s_1, A_1) \leq (s_2, A_2)$ iff $s_1 \subseteq s_2, A_1 \supseteq A_2$ and $s_2 \setminus s_1 \subseteq A_1$.

Let $G \subseteq \mathbb{M}_D$ be generic over ${\rm \bf V}$. The Mathias real $x_G$ is $\bigcup \{ s : \exists A \in D, (s, A) \in G \}$. Notice that in ${\rm \bf V}^{\mathbb{M}_D}$ we have $(x_G \subseteq^* B) \vee (x_G \subseteq^* \omega \setminus B)$ for every $B \in [\omega]^\omega$. We shall use this profound property while trying to create the monochromatic subsets in the theorems below. In \S 1 we use the original Mathias forcing, and in \S 2 we use the straightforward generalization of it for higher cardinals.

We thank the referees for many helpful comments, corrections and improvements.

\newpage 

\section{The countable case}
We prove, in this section, the consistency of $\binom{\mathfrak{c}}{\omega} \rightarrow \binom{\mathfrak{c}}{\omega}^{1,1}_n$, for every natural number $n$. The general pattern of the proof will be used also in the next section.

\begin{theorem}
\label{mt}
A positive realtion for $\mathfrak{c}$. \newline 
The strong relation $\binom{\mathfrak{c}}{\omega} \rightarrow \binom{\mathfrak{c}}{\omega}^{1,1}_n$ is consistent with ZFC, for every natural number $n$.
\end{theorem}

\par \noindent \emph{Proof}. \newline 
Choose an uncountable cardinal $\lambda$ so that $\lambda^{\aleph_0} = \lambda$ and $\cf(\lambda) \geq \aleph_2$. We define a finite support iteration $\langle \mathbb{P}_i, \name{\mathbb{Q}}_j : i \leq \omega_1, j < \omega_1 \rangle$ of $ccc$ forcing notions, such that $|\mathbb{P}_i| = \lambda$ for every $i \leq \omega_1$.

Let $\name{\mathbb{Q}}_0$ be (a name of) a forcing notion which adds $\lambda$ reals (e.g., Cohen forcing). For every $j < \omega_1$ let $\name{D}_j$ be a name of a nonprincipal ultrafilter on $\omega$. Let $\name{\mathbb{Q}}_{1+j}$ be a $ccc$ forcing notion which adds an infinite set $\name{A}_j \subseteq \omega$ such that $(\forall B \in \name{D}_j) (\name{A}_j \subseteq^* B \vee \name{A}_j \subseteq^* \omega \setminus B)$. The Mathias forcing $\mathbb{M}_{\name{D}_j}$ can serve. At the end, set $\mathbb{P} = \bigcup \{ \mathbb{P}_i : i < \omega_1 \}$.

Since every component satisfies the $ccc$, and we use finite support iteration, $\mathbb{P}$ is also a $ccc$ forcing notion and hence no cardinal is collapsed in ${\rm \bf V}^\mathbb{P}$. In addition, notice that $2^{\aleph_0} = \lambda$ after forcing with $\mathbb{P}$. Our goal is to prove that $\binom{\lambda}{\omega} \rightarrow \binom{\lambda}{\omega}^{1,1}_n$ in ${\rm \bf V}^\mathbb{P}$.

Let $\name{c}$ be a name of a function from $\lambda \times \omega$ into $n$. For every $\alpha < 2^{\aleph_0}$ we have a name (in ${\rm \bf V}$) to the restriction $\name{c} \upharpoonright (\{ \alpha \} \times \omega)$. $\mathbb{P}$ is $ccc$, hence the color of every pair of the form $(\alpha, n)$ is determined by an antichain which includes at most $\aleph_0$ conditions. Since we have to decide the color of $\aleph_0$-many pairs in $\name{c} \upharpoonright (\{ \alpha \} \times \omega)$, and the length of $\mathbb{P}$ is $\aleph_1$, we know that $\name{c} \upharpoonright (\{ \alpha \} \times \omega)$ is a name in $\mathbb{P}_{i(\alpha)}$ for some $i(\alpha) < \omega_1$.

For every $j < \omega_1$ let $\mathcal{U}_j$ be the set $\{ \alpha < \lambda : i(\alpha) \leq j \}$. Recall that $\aleph_1 < \aleph_2 \leq \cf(\lambda)$, hence for some $j < \omega_1$ we have $\mathcal{U}_j \in [\lambda]^\lambda$. Choose such $j$, and denote $\mathcal{U}_j$ by $\mathcal{U}$. We shall try to show that $\mathcal{U}$ can serve (after some shrinking) as the first coordinate in the monochromatic subset. 

Choose a generic subset $G \subseteq \mathbb{P}$, and denote $\name{A}_j[G]$ by $A$. For each $\alpha \in \mathcal{U}$ we know that $\name{c} \upharpoonright (\{ \alpha \} \times A)$ is constant, except a possible mistake over a finite subset of $A$. But this mistake can be amended.

For every $\alpha \in \mathcal{U}$ choose $k(\alpha) \in \omega$ and $m(\alpha) < n$ so that $(\forall \ell \in A)[\ell \geq k(\alpha) \Rightarrow \name{c}[G](\alpha, \ell) = m(\alpha)]$. $n$ is finite and $\cf(\lambda) > \aleph_0$, so one can fix some $k \in \omega$ and a color $m < n$ such that for some $\mathcal{U}_1 \in [\mathcal{U}]^\lambda$ we have $\alpha \in \mathcal{U}_1 \Rightarrow k(\alpha) = k \wedge m(\alpha) = m$.

Let $B$ be $A \setminus k$, so $B \in [\omega]^\omega$. By the fact that $\mathcal{U}_1 \subseteq \mathcal{U}$ we know that $c(\alpha, \ell) = m$ for every $\alpha \in \mathcal{U}_1$ and $\ell \in B$, yielding the positive relation $\binom{\mathfrak{c}}{\omega} \rightarrow \binom{\mathfrak{c}}{\omega}^{1,1}_n$, as required.

\hfill \qedref{mt}

\begin{remark}
\label{sstt}
Assume $\lambda$ is an uncountable regular cardinal. Denote by $\binom{\lambda}{\mu} \rightarrow_{\rm st} \binom{\lambda}{\mu}^{1,1}_n$ the assertion that for every coloring $c : \lambda \times \mu \rightarrow n$ there exists $B \in [\mu]^\mu$ and a stationary subset $\mathcal{U} \subseteq \lambda$ so that $c \upharpoonright \mathcal{U} \times B$ is constant. Our proof gives the consistency of $\binom{\mathfrak{c}}{\omega} \rightarrow_{\rm st} \binom{\mathfrak{c}}{\omega}^{1,1}_n$, when the continuum is regular.
\end{remark}

Let us turn to cardinals below the continuum. We deal with the relation $\binom{\theta}{\omega} \rightarrow \binom{\theta}{\omega}^{1,1}_2$, when $\theta < \mathfrak{s}$. We shall prove that this relation holds iff $\theta$ is of uncountable cofinality. 

Recall:

\begin{definition}
\label{ssplitting}
The splitting number. \newline 
Let $\mathcal{F} = \{ S_\alpha : \alpha < \kappa \}$ be a family of subsets of $\omega$, and $B \in [\omega]^\omega$.
\begin{enumerate}
\item [($\aleph$)] $\mathcal{F}$ splits $B$ if $|B \cap S_\alpha| = |B \cap \omega \setminus S_\alpha| = \aleph_0$ for some $\alpha < \kappa$
\item [($\beth$)] $\mathcal{F}$ is a splitting family if $\mathcal{F}$ splits $B$ for every $B \in [\omega]^\omega$
\item [($\gimel$)] $\mathfrak{s} = {\rm min} \{ |\mathcal{F}| : \mathcal{F} {\rm \ is \ a  \ splitting \ family} \}$
\end{enumerate}
\end{definition}

\begin{claim}
\label{ppss}
The polarized relation below $\mathfrak{s}$. \newline 
Assume $\omega \leq \theta < \mathfrak{s}$. \newline 
\Then\ $\binom{\theta}{\omega} \rightarrow \binom{\theta}{\omega}^{1,1}_2$ iff
$\cf(\theta) > \aleph_0$.
\end{claim}

\par \noindent \emph{Proof}. \newline 
Suppose first that $\cf(\theta) > \aleph_0$. Let $c : \theta \times \omega \rightarrow 2$ be a coloring. Define $S_\alpha = \{ n \in \omega : c(\alpha, n) = 0 \}$ for every $\alpha < \theta$, and $\mathcal{F} = \mathcal{F}_c = \{ S_\alpha : \alpha < \theta \}$. Since $\theta < \mathfrak{s}$, we know that $\mathcal{F}$ is not a splitting family.

Choose an evidence, i.e., $B \in [\omega]^\omega$ which is not splitted by $\mathcal{F}$. It means that $(B \subseteq^* S_\alpha) \vee (B \subseteq^* \omega \setminus S_\alpha)$ for every $\alpha < \theta$. At least one of this two opstions occurs $\theta$-many times, and since all we need for the first coordinate (in the monochromatic subset) is its cardinality, we shall assume (without loss of generality) that $B \subseteq^* S_\alpha$ for every $\alpha < \theta$.

For every $\alpha < \theta$ there exists a finite set $t_\alpha \subset \omega$ such that $B \setminus t_\alpha \subseteq S_\alpha$. There are countably-many $t_\alpha$-s, and $\cf(\theta) > \aleph_0$, so for some $t \in [\omega]^{< \omega}$ and $H_0 \in [\theta]^\theta$ we have $\alpha \in H_0 \Rightarrow B \setminus t \subseteq S_\alpha$. Set $H_1 = B \setminus t$, and verify that $c \upharpoonright H_0 \times H_1 \equiv 0$, so we are done.

Now suppose $\cf(\theta) = \aleph_0$, and choose an increasing sequence $\langle \theta_n : n \in \omega \rangle$ which tends to $\theta$. For every $\alpha < \theta$, let $\ell(\alpha)$ be the first natural number $n$ such that $\theta_n \leq \alpha < \theta_{n + 1}$. Define $c(\alpha, n) = 0 \Leftrightarrow \ell(\alpha) \geq n$. We claim that $c$ is an evidence to the negative assertion $\binom{\theta}{\omega} \nrightarrow \binom{\theta}{\omega}^{1,1}_2$.

Indeed, assume $H \in [\theta]^\theta$ and $B \in [\omega]^\omega$. If $c \upharpoonright H \times B \equiv 1$ then $\ell(\alpha) < n$ for every $\alpha \in H$ and $n \in B$. Choose some specific $n \in B$. Since $H$ is unbounded in $\theta$, one can pick large enough $\alpha' \in H$ such that $\ell(\alpha') \geq n$. Consequently, $c(\alpha', n) = 0$, a contradiction. On the other hand, if $c \upharpoonright H \times B \equiv 0$ then $\ell(\alpha) \geq n$ for every $\alpha \in H$ and $n \in B$. Choose some specific $\alpha \in H$. Since $B$ is unbounded in $\omega$, one can pick large enough $n' \in B$ such that $\ell(\alpha) < n'$. Consequently, $c(\alpha, n') = 1$, a contradiction. So the proof is complete.

\hfill \qedref{ppss}

We indicate that for $\mathfrak{s}$ itself we believe that the relation $\binom{\mathfrak{s}}{\omega} \rightarrow \binom{\mathfrak{s}}{\omega}^{1,1}_2$ is independent of ZFC. We hope to shed light on this issue in a subsequent work.

\newpage 

\section{The supercompact case}

In this section we prove the consistency of $\binom{2^\mu}{\mu} \rightarrow \binom{2^\mu}{\mu}^{1,1}_2$ for every supercompact $\mu$. We shall use \cite{MR0472529} for making $\mu$ indestructible (in fact, all we need is the measurability of $\mu$ at every stage of the iteration), and \cite{MR0505492} for preserving the property of $\lambda$-cc along the iteration. We shall use a generalization of the Mathias forcing, so we need the following:

\begin{definition}
\label{ggeneralmat}
The generalized Mathias forcing. \newline 
Let $\mu$ be a supercompact (or even just measurable) cardinal, and $D$ a nonprincipal $\mu$-complete ultrafilter on $\mu$. The forcing notion $\mathbb{M}_D^\mu$ consists of pairs $(a, A)$ such that $a \in [\mu]^{< \mu}, A \in D$. For the order, $(a_1, A_1) \leq (a_2, A_2)$ iff $a_1 \subseteq a_2, A_1 \supseteq A_2$ and $a_2 \setminus a_1 \subseteq A_1$.
\end{definition}

If $\mathbb{M}_D^\mu$ is a $\mu$-Mathias forcing, then for defining the Mathias $\mu$-real we take a generic $G \subseteq \mathbb{M}_D^\mu$, and define $x_G = \bigcup \{ a : (\exists A \in D)((a, A) \in G) \}$. As in the original Mathias forcing, $x_G$ is endowed with the property $x_G \subseteq^* A \vee x_G \subseteq^* \mu \setminus A$ for every $A \in [\mu]^\mu$.

\begin{theorem}
\label{ssupcomp}
A positive relation for $2^\mu$. \newline 
The strong relation $\binom{2^\mu}{\mu} \rightarrow \binom{2^\mu}{\mu}^{1,1}_\theta$ is consistent with ZFC, for every $\theta < \mu$.
\end{theorem}

\par \noindent \emph{Proof}. \newline 
Let $\mu$ be a supercompact cardinal. Starting with Laver's forcing, we may assume that $\mu$ is Laver-indestructible. Choose any $\lambda$ so that $\cf(\lambda) = \mu^{++}$, and $\lambda^\mu = \lambda$. Let $\langle \mathbb{P}_\alpha, \name{\mathbb{Q}}_\beta : \alpha \leq \mu^+, \beta < \mu^+ \rangle$ be an iteration with $(< \mu)$-support. We start with $\mathbb{P}_0 = \{ \emptyset \}$ and $\name{\mathbb{Q}}_0$ a name in $\mathbb{P}_0$ of a forcing which increases $2^\mu$ to $\lambda$ (e.g., Cohen forcing).

For every $\beta < \mu^+$ we choose $\name{D}_\beta$, a name of a nonprincipal $\mu$-complete ultrafilter on $\mu$. Let $\name{\mathbb{Q}}_{1 + \beta}$ be (a name of) the generalized Mathias forcing $\mathbb{M}^\mu_{\name{D}_\beta}$. Notice that $\mathbb{M}^\mu_{\name{D}_\beta}$ is $\mu$-closed (since $\name{D}_\beta$ is a $\mu$-complete ultrafilter), so $\mu$ remains supercompact and hence measurable along the iteration.

For every $0 < \beta < \mu^+$, choose a generic set $\name{G}_{1 + \beta} \subseteq \name{\mathbb{Q}}_{1 + \beta}$, and let $\name{A}_{1 + \beta}$ be the Mathias $\mu$-real associated with it. We shall work in ${\rm \bf V}^{\mathbb{P}}$, when $\mathbb{P} = \bigcup \{ \mathbb{P}_\alpha : \alpha < \mu^+ \}$, aiming to show the positive relation $\binom{2^\mu}{\mu} \rightarrow \binom{2^\mu}{\mu}^{1,1}_\theta$. 

First, let us indicate that $\mathbb{P}$ satisfies the $\mu^+$-cc. It follows from \cite{MR0505492}, upon noticing that each component satisfies a strong form of the $\mu^+$-cc as required there. Second, $\mathbb{P}$ is $\mu$-complete, since each component is $\mu$-complete. Consequently, no cardinal is collapsed and no cofinality is changed by $\mathbb{P}$. Moreover, $\name{\mathbb{Q}}_0$ blows $2^\mu$ to $\lambda$, and the completeness of the other forcing notions ensures that ${\rm \bf V}^{\mathbb{P}} \models 2^\mu = \lambda$. We shall prove that $\binom{\lambda}{\mu} \rightarrow \binom{\lambda}{\mu}^{1,1}_\theta$ in ${\rm \bf V}^{\mathbb{P}}$.

Assume that $\theta < \mu$ is fixed, and $\name{c}$ is a name of a coloring function from $\lambda \times \mu$ into $\theta$. We denote $\name{c} \upharpoonright ( \{ \alpha \} \times \mu )$ by $\name{c}_\alpha$, for every $\alpha < \lambda$, and we claim that $\name{c}_\alpha \in \mathbb{P}_{\xi(\alpha)}$ for some $\xi(\alpha) < \mu^+$. For this, notice that $\{ \alpha \} \times \mu$ consists of $\mu$ pairs, and for the color of each pair we have at most $\mu$ conditions which give different values, since $\mathbb{P}$ is $\mu^+$-cc. But $\mathbb{P}$ is of length $\mu^+$, so $\name{c}_\alpha$ appears at some early stage $\mathbb{P}_{\xi(\alpha)}$.

For every $\beta < \mu^+$, set $\mathcal{U}_\beta = \{ \alpha < \lambda : \xi(\alpha) \leq \beta \}$. Since $\cf(\lambda) = \mu^{++} > \mu^+$, one can pick an ordinal $\beta < \mu^+$ so that $|\mathcal{U}_\beta| = \lambda$. Let $\mathcal{U}$ be $\mathcal{U}_\beta$, and let $G \subseteq \mathbb{P}$ be generic over ${\rm \bf V}$. Denote $\name{A}_\beta[G]$ by $A$.

For every $\alpha \in \mathcal{U}$, $\name{c}_\alpha$ is constant on $A$, except a small (i.e., of cardinality less than $\mu$) subset of $A$. For each $\alpha \in \mathcal{U}$ choose $\zeta(\alpha) \in \mu$ and $\theta(\alpha) \in \theta$ such that $\zeta \geq \zeta(\alpha) \wedge \zeta \in A \Rightarrow \name{c}[G](\alpha, \zeta) = \theta(\alpha)$. By the assumptions on $\theta, \mu$ and the cofinality of $\lambda$, there is $\mathcal{U}_1 \in [\mathcal{U}]^\lambda$ and $\zeta_* \in \mu, \theta_* \in \theta$ such that $\zeta(\alpha) \equiv \zeta_*$ and $\theta(\alpha) \equiv \theta_*$ for every $\alpha \in \mathcal{U}_1$.

Set $B = A \setminus \zeta_*$, and notice that $|B| = \mu$. By the above considerations, for every $\alpha \in \mathcal{U}_1$ and every $\zeta \in B$ we have $c(\alpha, \zeta) = \theta_*$. Hence $\binom{\lambda}{\mu} \rightarrow \binom{\lambda}{\mu}^{1,1}_\theta$, and the proof is complete.

\hfill \qedref{ssupcomp}

\begin{remark}
\label{refff}
Notice that the relation $\binom{\lambda}{\mu} \rightarrow \binom{\lambda}{\mu}^{1,1}_2$ is completely determined under GCH, as follows from \cite{williams} theorem 4.2.8.
\end{remark}

\newpage 

\bibliographystyle{amsplain}
\bibliography{arlist}

\providecommand{\bysame}{\leavevmode\hbox to3em{\hrulefill}\thinspace}
\providecommand{\MR}{\relax\ifhmode\unskip\space\fi MR }
\providecommand{\MRhref}[2]{%
  \href{http://www.ams.org/mathscinet-getitem?mr=#1}{#2}
}
\providecommand{\href}[2]{#2}
\begin{thebibliography}{1}

\bibitem{MR0202613}
P.~Erd{\H{o}}s, A.~Hajnal, and R.~Rado, \emph{Partition relations for cardinal
  numbers}, Acta Math. Acad. Sci. Hungar. \textbf{16} (1965), 93--196.
  \MR{MR0202613 (34 \#2475)}

\bibitem{MR0081864}
P.~Erd{\"o}s and R.~Rado, \emph{A partition calculus in set theory}, Bull.
  Amer. Math. Soc. \textbf{62} (1956), 427--489. \MR{MR0081864 (18,458a)}

\bibitem{GaSh949}
Shimon Garti and Saharon Shelah, \emph{A strong polarized relation}, in
  preperation.

\bibitem{MR0371652}
R.~Laver, \emph{Partition relations for uncountable cardinals {$\leq 2^{\aleph
  _{0}}$}}, Infinite and finite sets ({C}olloq., {K}eszthely, 1973; dedicated
  to {P}. {E}rd{\H o}s on his 60th birthday), {V}ol. {II}, North-Holland,
  Amsterdam, 1975, pp.~1029--1042. Colloq. Math. Soc. Jan\'os Bolyai, Vol. 10.
  \MR{MR0371652 (51 \#7870)}

\bibitem{MR0472529}
Richard Laver, \emph{Making the supercompactness of {$\kappa $} indestructible
  under {$\kappa $}-directed closed forcing}, Israel J. Math. \textbf{29}
  (1978), no.~4, 385--388. \MR{MR0472529 (57 \#12226)}

\bibitem{MR0505492}
S.~Shelah, \emph{A weak generalization of {MA} to higher cardinals}, Israel J.
  Math. \textbf{30} (1978), no.~4, 297--306. \MR{MR0505492 (58 \#21606)}

\bibitem{williams}
Neil~H. Williams, \emph{Combinatorial set theory}, studies in logic and the
  foundations of mathematics, vol.~91, North-Holland publishing company,
  Amsterdam, New York, Oxford, 1977.

\end{thebibliography}

\end{document}